\documentclass{article}
\usepackage{amsmath}
\usepackage{longtable}
\usepackage{url}
\usepackage[margin=10truemm]{geometry}
\title{\huge  On\ $A^4 + hB^4 = C^4 + hD^4$}
\author{Seiji Tomita }
\date{}
\begin{document}
\maketitle
\begin{abstract}
In this paper, we proved that $A^4 + hB^4 = C^4 + hD^4$ always has the integral solutions for $h < 20000.$
\end{abstract}
\vskip\baselineskip

\centerline{\huge1. Introduction}
\vskip\baselineskip

E. Grigorief\cite{d} noted $19^4+5\cdot281^4=417^4+5\cdot117^4$ and he found an infinitude of solutions when $h=2.$

A.S.Werebrusow\cite{d} gave $139^4+2\cdot34^4=61^4+2\cdot116^4.$

A. Gerardin\cite{d} gave the solution $(a,b,c,d)=(2p^2,p-1,2p,p+1)$ for $h=2p^3(p^2-1)$ and gave solutions for

 $26$ numerical values of $h.$

In 1997, Ajai Choudhry\cite{b} gave the integer solutions of $A^4+hB^4=C^4+hD^4$ with $h \leq 101.$

Furthermore, he proved that using given non trivial solution, other non trivial solution can be found.

In 2013, this author extended the solution table for  $h < 1000$ except $h=967$ and gave several parametric solutions

 when $h$ is given by polynomials.

The missing solution for $h = 967$ was supplied by Andrew Bremner\cite{a}.

In 2016, Ajai Chodhry\cite{c} showed the parametric solutions when $h$ is given by polynomials of degrees $2, 3$, and $4.$

In 2017, Jaroslaw Wroblewski\cite{h} kindly sent me many solutions for  $h < 20000.$

In 2017, Noam Elkies\cite{e} found the solutions for $h=9719$ and $16329.$

In 2017, this author extended the solution table for  $h < 20000.$

\vskip\baselineskip

We searched the solutions by several methods as follows. 
 
The case $A<50000:$ Brute force  

The case $A<100000:$ Sorted sum  

The case $A<1000000:$ Quartic  

The case $A>1000000:$ Elliptic curve  
\vskip\baselineskip

Now, we found the all solutions for $h < 20000.$  

This contains all presently known solutions.

We show the numerical solutions only for $h<1000$, other solutions are obtained from the Appendix.
\vskip\baselineskip

\newpage

\centerline{\huge2. Brute force}
\vskip3\baselineskip

We consider to find the integer solutions of $A^4+hB^4=C^4+hD^4$.

Since $A^4-C^4=h(D^4-B^4),$ then $A^4-C^4 \equiv 0 \pmod h.$

Let $A > B,C$ and $A<50000.$

If $h$ has a prime factor $p$ such that $p \equiv 3 \pmod 4,$ then $A^2+C^2 \not\equiv  0 \pmod p.$

Since $p$ cannot divide $A^2+C^2,$ so divides $A^2-C^2.$

Hence  $p$ divides $A-C$ or $A+C.$ 

First we consider the case of $p$ divides $A-C$.

We set the search range such as $A_{min}<A<A_{max}$ and $C_{min}<C<C_{max}$.

Hence let $C$ range from $C_{min}$ and take $A$ such that $A=C+[A_{min}-C)/h]h+hm, m=0,1,2,\cdots.$

\vskip\baselineskip
We tested whether $f=\cfrac{(A^4-C^4 +hB^4)}{h}$ could be written in the form $x^4.$ 
\vskip\baselineskip

If $f$ is written in the form $x^4,$ we can obtain the integer solutions.
\vskip\baselineskip

We used the useful relation as follows to reduce the calculation time.

If $f = x^4$ then $f \equiv 0,1 \pmod {16}.$
\vskip\baselineskip

Let $h=4117=23\cdot179 ,\ A_{min}=10000,\ B=2303,\ C=2263,\ p=4117,\ m=0,1,2,\cdots$ then we obtain the sequence,
\vskip\baselineskip

$C+[A_{min}-C)/h]h+hm=2263+4117,\ 2263+4117+4117\cdot1,\ 2263+4117+4117\cdot2, \cdots$
\vskip\baselineskip

Thus we know  $A=2263+4117+4117=10497$,\ $f=\cfrac{(A^4-C^4 +hB^4)}{h}=2361^4$, then we obtain
\vskip\baselineskip

$(A,B,C,D)=(10497,\ 2303,\ 2263,\ 2361).$
\vskip\baselineskip

Next we consider the case of $p$ divides $A+C$.

We searched the solutions same as for the case of $p$ divides $A-C$.

\vskip3\baselineskip

\centerline{\huge3. Sorted sum}
\vskip2\baselineskip

We consider to find the integer solutions of $A^4+hB^4=C^4+hD^4$ with $50000 < A \leq 100000$.
\vskip\baselineskip

Sort the set {$(A^4+hB^4, A, B) : A, B < 100000$} according to $A^4+hB^4$ and find the repeated entries.
\vskip\baselineskip

Let $p,q$ are prime number such that $p,q \equiv 3 \pmod 4$  

For the given $rp$ where $0 \leq rp<p$, generate the set {($A^4+hB^4, A, B$) } where $A^4+hB^4 \equiv rp \pmod p$.
\vskip\baselineskip

For the given $rq$ where $0 \leq rq<q$, choose the set {($A^4+hB^4, A, B$) } $where A^4+hB^4 \equiv rq \pmod q$.
\vskip\baselineskip

Next, sort {($A^4+hB^4, A, B$) } according to $A^4+hB^4$ and find the repeated entries.
\vskip\baselineskip

Example:  $A^4 + 2518B^4 = C^4 + 2518D^4$
\vskip\baselineskip

Let $p=331,q=347, rp=304, rq=35$, we obtain one of the candidates
\vskip\baselineskip

$58948^4+2518\cdot5687^4 \equiv 61916^4+2518\cdot1481^4 \equiv 304 \pmod {331}$
\vskip\baselineskip

$58948^4+2518\cdot5687^4 \equiv 61916^4+2518\cdot1481^4 \equiv 35 \pmod {347}$
\vskip\baselineskip

Since this candidate certainly satisfies the equation $A^4+hB^4=C^4+hD^4$, then we obtain
\vskip\baselineskip

$(A,B,C,D)=(58948,\ 5687,\ 61916,\ 1481)$.
\newpage

\centerline{\huge4. Quartic and Elliptic curve}
\vskip3\baselineskip
Equation can be reduced to Quartic and elliptic curve.
\vskip\baselineskip

Let $A = px+a, B = qx-b, C = px-a, D = qx+b,$ then we obtain

$$A^4 + hB^4 - C^4 - hD^4 = (8ap^3-8hbq^3)x^2+(-8hb^3q+8a^3p)$$

Since $x$ must be rational number, then discriminant for $x$  must be square number.

$$v^2 = -a^4p^4+hb^3qap^3+a^3hbq^3p-h^2b^4q^4$$

Let $U=\frac{p}{q}, V=\frac{v}{q^2}$ then we obtain Quartic

$$V^2 = -a^4U^4+hb^3aU^3+a^3hbU-h^2b^4$$

Hence we search the rational solutions $(U,V)$ for given $(a,b,h).$
\vskip\baselineskip

We can find the rational solution of quartic using \textbf{PARI-GP\cite{f}: hyperellratpoints(Michael Stoll)}.
\vskip\baselineskip

We searched the rational solution of quartic where $100000<A<1000000.$
\vskip\baselineskip

Let $h=2572$ , $a=15$, and $b=14$ then 
\vskip\baselineskip

We found a rational solution $(U,V)=(\cfrac{9002}{267}, \cfrac{45966408488}{23763})$, then we obtain
\vskip\baselineskip

$(A,B,C,D)=(799298,\ 61171,\ 623018,\ 103357).$
\vskip\baselineskip

If we could not find the rational solution, we used elliptic curve method as follows.

Quartic can be reduced to elliptic curve below.

$$Y^2 = X^3-3h^2b^4a^4X-b^2h^2a^2(a^8+b^8h^2)$$
with
$$U = \frac{-1}{2}\frac{2bha^6+2b^3hX}{-ab^6h^2-a^3X}$$
$$V = \frac{bh(a^8-b^8h^2)Y}{ab^{12}h^4+2Xa^3h^2b^6+X^2a^5}$$

We can find the rational solution of elliptic curve using \textbf{SAGE\cite{g}: gens}.

Let $h=9069$ , $a=3$, and $b=1$ then 

$$Y^2 = X^3-19985962923X-60885623843910378$$

The elliptic curve has rank $2$ and generators = 
\begin{align*}
\Bigl(&[\frac{11633949063}{14161}, \frac{1164093129464040}{1685159}],\\ 
 &[\frac{4587889797054}{6723649}, \frac{-8597517313555330650}{17434421857}]\Bigr)
\end{align*}

\vskip\baselineskip

Take $X=\cfrac{11633949063}{14161}$ then $U = \cfrac{74903894}{2701177}.$
\vskip\baselineskip

Thus $p = 74903894$ and $q = 2701177$ then we obtain $x = \cfrac{119}{825696345}.$
\vskip\baselineskip

Hence $(A,B,C,D)=(11390652421,\ 504256282,\ 6436474351,\ 1147136408).$
\newpage

\centerline{\huge5. Parametric solutions}
\vskip\baselineskip

When $h = a^4 + b^4$, equation has obvious solution $(b^2)^4 + ha^4 = (a^2)^4 + hb^4.$

Gerardin gave a parametric solution below.

$(2n^2)^4 + h(n-1)^4 = (2n)^4 + h(n+1)^4$ where $h=2n^3(n^2-1).$

In 2016, Ajai Choudhry found the parametric solutions when $h$ is  given by polynomials of degrees 2, 3, and 4.
\begin{table}[hbtp]
\centering
  \caption{Solutions of $A^4 + hB^4 = c^4 + hD^4$}
  \label{table:data_type}
  \begin{tabular}{ccccc}
    \hline
   h & A  & B  &  C  & D   \\
    \hline \hline

$p^2 + 2$  & $p^3 + 2p + 1$  & $p^2- p + 1$ & $p^3 + 2p - 1$ & $p^2 + p + 1$ \\
$p(p^2 + 4)$ & $p -2$ & $2$ & $p + 2$ &  $0$ \\
$8p(p^2 + 1)$ & $p - 1$ & $1$ & $p + 1$ & $0$ \\
$p^4- 1$ & $p$ & $0$ & $1$ & $1$ \\
$2p^4- 2$ & $p^2 + 2p - 1$ & $p - 1$ & $p^2- 2p - 1$ & $p + 1$ \\
$p^4 + 3p^2 + 1$ & $p^2 + p + 1$ & $p - 1$ & $p^2- p + 1$  & $p + 1$ \\

    \hline
    $p$ is arbitrary.
  \end{tabular}
\end{table}

In 2016, this author found the parametric solutions when h is  given by polynomials of degrees 4, 6, and 7.

\begin{table}[hbtp]
\centering
  \caption{Solutions of $A^4 + hB^4 = c^4 + hD^4$}
  \label{table:data_type}
  \begin{tabular}{ccccc}
    \hline
   h & A  & B  &  C  & D   \\
    \hline \hline

 $ (n^2+4)(n^2+2)  $ &                  $n^2+n+2  $ &  $n-1  $ &      $n^2-n+2  $ &  $n+1$  \\

 $ 2n^4+12n^2+2  $ &                    $n^2+2n+1  $ & $n-1  $ &     $ n^2-2n+1  $ & $n+1   $ \\

  $2(n^2+9)(n^2+3)  $ &                 $n^2+2n+3  $ & $n-1  $ &     $ n^2-2n+3  $ & $n+1  $ \\

  $3(n^2+4)(n^2-2)  $ &                 $n^2+3n-2  $ & $n-1  $ &     $ n^2-3n-2  $ & $n+1  $ \\

  $3n^4+33n^2+3  $ &                    $n^2+3n+1  $ & $n-1  $ &      $n^2-3n+1  $ & $n+1  $ \\

  $n^6+2n^4+n^2+1  $ &                  $n^3+n+1  $ &  $n-1  $ &      $n^3+n-1  $ &  $n+1  $ \\

  $2n^6+4n^4+2n^2+8  $ &                $n^3+n+2  $ &  $n-1  $ &      $n^3+n-2  $ &  $n+1  $ \\

  $3n^6+6n^4+3n^2+27  $ &               $n^3+n+3  $ &  $n-1  $ &     $ n^3+n-3  $ &  $n+1  $ \\

  $(n^2+2)(n^4+3n^2+1)  $ &             $n^3+2n+1  $ & $n-1  $ &     $ n^3+2n-1  $ & $n+1  $ \\

  $(2n+1)(n^6+2n^4+5n^2+4n+1)  $ &      $n^3+3n+1  $ & $n-1  $ &      $n^3-n-1  $ &  $n+1  $ \\

  $2(n+1)(n^6+2n^4+5n^2+8n+4)  $ &      $n^3+3n+2  $ & $n-1  $ &      $n^3-n-2  $ &  $n+1  $ \\

  $2(n+1)(n^2-n+2)(n^4+n^2+4n+4)  $ &   $n^3+3n+2  $ & $n-1  $ &     $ n^3-n+2  $ &  $n+1  $ \\

  $2(n^2+4)(n^2+3)(n^2+1)  $ &         $ n^3+3n+2  $ & $n-1  $ &     $ n^3+3n-2  $ & $n+1  $ \\

  $(n^2+3)(n^2+1)^2  $ &               $ n^3+3n+2  $ & $n-2  $ &     $ n^3+3n-2  $ & $n+2  $ \\

  $(2n+3)(n^6+2n^4+5n^2+12n+9)  $ &    $ n^3+3+3n  $ & $n-1  $ &     $ n^3-n-3  $ & $ n+1  $ \\

    \hline
    
  $n$ is arbitrary.
  \end{tabular}
\end{table}

We show two more solutions.
\vskip\baselineskip

$(a): h = 2(m^2-n^2)^3mn$

$$v^2 = -a^4p^4+hb^3qap^3+a^3hbq^3p-h^2b^4q^4$$

Substitute $h=\frac{p^3a}{q}$ and  $b=1$ to quartic, then

$$v^2 = -p^4a^2(p-a)(p+a)(q-1)(q+1)$$

Let $a = m^2+n^2, p = m^2-n^2, q = \frac{m^2+n^2}{2mn},$ then $h = 2(m^2-n^2)^3mn$ and  $x = \frac{4m^2n^2}{-2m^2n^2+n^4+m^4}.$
\vskip\baselineskip

Hence we obtain a parametric solution as follows.
\vskip\baselineskip

$A = (n-m)(n+m)(n^4-4m^2n^2-m^4)$
\vskip\baselineskip

$B = 2m^3n+2n^3m+2m^2n^2-n^4-m^4$
\vskip\baselineskip

$C = (n-m)(n+m)(n^4+4m^2n^2-m^4)$
\vskip\baselineskip

$D = 2m^3n+2n^3m-2m^2n^2+n^4+m^4$

$m,n$ are arbitrary.

\vskip\baselineskip

$(b): h = 8(m^2-n^2)m^3n^3$
\vskip\baselineskip

Substitute $h=\frac{pa}{q}$ and  $b=1$ to quartic, then

$$v^2 = -a^2p^2(a-1)(a+1)(p-q)(p+q)$$

Let $a = \frac{m^2+n^2}{2mn}, q = m^2+n^2, p = m^2-n^2,$ then $h = 8(m^2-n^2)m^3n^3$ and  $x = \frac{-(-m^2+n^2)}{4(m^2n^2)} .$
\vskip\baselineskip

Hence we obtain a parametric solution as follows.
\vskip\baselineskip

$A = 2mn(m^4-2m^2n^2+n^4+2m^3n+2mn^3)$
\vskip\baselineskip

$B = m^4-n^4-4m^2n^2$
\vskip\baselineskip

$C = 2mn(m^4-2m^2n^2+n^4-2m^3n-2mn^3)mn$
\vskip\baselineskip

$D = m^4-n^4+4m^2n^2$

$m,n$ are arbitrary.
\newpage

\centerline{\huge6. Numerical solutions}
\vskip3\baselineskip

\begin{align*}
&\text{ Smallest: Min}(A^4 + hB^4)\\  
&\text{ Search range:} h<1000, (A,B,C,D)<100000  \\
&\text{ When h is fourth power, it is excluded from searching} \\
&\text{ ???: May not be smallest solution}
\end{align*}

\begin{longtable}[c]{ccccc}
\caption{\Large{Solution table}}
\label{longtablesample} \\
\hline
h & A & B & C & D \\
\hline \endhead

2 & 139 & 34 & 61 & 116  \\ 
3 & 4 & 1 & 2 & 3  \\ 
4 & 9 & 4 & 7 & 6  \\ 
5 & 3 & 0 & 1 & 2  \\ 
6 & 13 & 3 & 11 & 7  \\ 
7 & 13 & 2 & 1 & 8  \\ 
8 & 232 & 61 & 68 & 139  \\ 
9 & 29 & 9 & 11 & 17  \\ 
10 & 7 & 2 & 1 & 4  \\ 
11 & 21 & 12 & 1 & 14  \\ 
12 & 31 & 9 & 1 & 17  \\ 
13 & 6 & 1 & 4 & 3  \\ 
14 & 11 & 4 & 3 & 6  \\ 
15 & 16 & 8 & 13 & 9  \\ 
17 & 4 & 1 & 1 & 2  \\ 
18 & 33 & 2 & 9 & 16  \\ 
19 & 98 & 13 & 16 & 47  \\ 
20 & 259 & 77 & 13 & 127  \\ 
21 & 41 & 8 & 29 & 18  \\ 
22 & 79 & 59 & 31 & 61  \\ 
23 & 29 & 14 & 17 & 16  \\ 
24 & 82 & 3 & 22 & 37  \\ 
25 & 149 & 68 & 43 & 80  \\ 
26 & 77 & 2 & 53 & 32  \\ 
27 & 9 & 2 & 3 & 4  \\ 
28 & 136 & 27 & 108 & 53  \\ 
29 & 44 & 7 & 14 & 19  \\ 
30 & 7 & 1 & 1 & 3  \\ 
31 & 18 & 8 & 13 & 9  \\ 
32 & 139 & 17 & 61 & 58  \\ 
33 & 7 & 2 & 4 & 3  \\ 
34 & 5 & 0 & 3 & 2  \\ 
35 & 8 & 1 & 6 & 3  \\ 
36 & 17 & 3 & 1 & 7  \\ 
37 & 26 & 7 & 8 & 11  \\ 
38 & 72 & 3 & 4 & 29  \\ 
39 & 5 & 0 & 1 & 2  \\ 
40 & 33 & 1 & 31 & 9  \\ 
41 & 48 & 9 & 22 & 19  \\ 
42 & 83 & 2 & 43 & 32  \\ 
43 & 367 & 257 & 23 & 263  \\ 
44 & 46 & 13 & 2 & 19  \\ 
45 & 19 & 4 & 17 & 6  \\ 
46 & 117 & 59 & 67 & 63  \\ 
47 & 24 & 6 & 23 & 7  \\ 
48 & 8 & 1 & 4 & 3  \\ 
49 & 44 & 17 & 5 & 20  \\ 
50 & 209 & 4 & 87 & 78  \\ 
51 & 8 & 1 & 2 & 3  \\ 
52 & 9 & 1 & 7 & 3  \\ 
53 & 31 & 2 & 29 & 8  \\ 
54 & 111 & 11 & 9 & 41  \\ 
55 & 14 & 1 & 8 & 5  \\ 
56 & 22 & 7 & 6 & 9  \\ 
57 & 49 & 8 & 11 & 18  \\ 
58 & 13 & 2 & 11 & 4  \\ 
59 & 654 & 93 & 644 & 127  \\ 
60 & 196 & 13 & 128 & 67  \\ 
61 & 11 & 2 & 1 & 4  \\ 
62 & 233 & 1 & 47 & 83  \\ 
63 & 8 & 2 & 1 & 3  \\ 
64 & 9 & 2 & 7 & 3  \\ 
65 & 3 & 0 & 2 & 1  \\ 
66 & 37 & 12 & 29 & 14  \\ 
67 & 179 & 49 & 89 & 67  \\ 
68 & 23 & 2 & 7 & 8  \\ 
69 & 32 & 1 & 14 & 11  \\ 
70 & 509 & 96 & 387 & 164  \\ 
71 & 236 & 67 & 94 & 89  \\ 
72 & 32 & 3 & 4 & 11  \\ 
73 & 38 & 1 & 4 & 13  \\ 
74 & 17 & 4 & 9 & 6  \\ 
75 & 19 & 1 & 17 & 5  \\ 
76 & 202 & 43 & 178 & 59  \\ 
77 & 248 & 109 & 214 & 113  \\ 
78 & 17 & 4 & 7 & 6  \\ 
79 & 41 & 5 & 38 & 10  \\ 
80 & 44 & 17 & 8 & 19  \\ 
82 & 9 & 1 & 1 & 3  \\ 
83 & 748 & 73 & 746 & 91  \\ 
84 & 31 & 0 & 17 & 10  \\ 
85 & 64 & 21 & 18 & 25  \\ 
86 & 79 & 8 & 7 & 26  \\ 
87 & 113 & 3 & 7 & 37  \\ 
88 & 158 & 51 & 18 & 61  \\ 
89 & 16 & 1 & 10 & 5  \\ 
90 & 27 & 8 & 3 & 10  \\ 
91 & 34 & 1 & 8 & 11  \\ 
92 & 1319 & 29 & 1303 & 199  \\ 
93 & 247 & 2 & 187 & 72  \\ 
94 & 111 & 16 & 17 & 36  \\ 
95 & 18 & 6 & 1 & 7  \\ 
96 & 26 & 3 & 22 & 7  \\ 
97 & 9 & 2 & 4 & 3  \\ 
98 & 119 & 47 & 63 & 51  \\ 
99 & 129 & 13 & 3 & 41  \\ 
100 & 53 & 9 & 21 & 17  \\ 
101 & 11 & 1 & 9 & 3  \\ 
102 & 103 & 19 & 47 & 33  \\ 
103 & 107 & 15 & 99 & 25  \\ 
104 & 146 & 53 & 62 & 59  \\ 
105 & 16 & 1 & 2 & 5  \\ 
106 & 37 & 12 & 3 & 14  \\ 
107 & 396 & 153 & 32 & 167  \\ 
108 & 93 & 11 & 3 & 29  \\ 
109 & 13 & 2 & 7 & 4  \\ 
110 & 303 & 36 & 71 & 94  \\ 
111 & 7 & 0 & 5 & 2  \\ 
112 & 13 & 1 & 1 & 4  \\ 
113 & 8 & 1 & 7 & 2  \\ 
114 & 31 & 9 & 7 & 11  \\ 
115 & 62 & 7 & 16 & 19  \\ 
116 & 88 & 11 & 12 & 27  \\ 
117 & 42 & 7 & 36 & 11  \\ 
118 & 269 & 8 & 203 & 74  \\ 
119 & 13 & 2 & 1 & 4  \\ 
120 & 16 & 3 & 4 & 5  \\ 
121 & 499 & 1 & 469 & 103  \\ 
122 & 447 & 29 & 207 & 133  \\ 
123 & 19 & 4 & 7 & 6  \\ 
124 & 78 & 31 & 46 & 33  \\ 
125 & 10 & 1 & 0 & 3  \\ 
126 & 187 & 57 & 61 & 67  \\ 
127 & 168 & 31 & 86 & 51  \\ 
128 & 257 & 11 & 223 & 62  \\ 
129 & 152 & 11 & 149 & 24  \\ 
130 & 17 & 1 & 7 & 5  \\ 
131 & 139 & 17 & 123 & 33  \\ 
132 & 28 & 7 & 16 & 9  \\ 
133 & 73 & 12 & 3 & 22  \\ 
134 & 296 & 9 & 28 & 87  \\ 
135 & 18 & 1 & 12 & 5  \\ 
136 & 38 & 3 & 18 & 11  \\ 
137 & 91 & 29 & 79 & 31  \\ 
138 & 556 & 61 & 88 & 163  \\ 
139 & 301 & 63 & 23 & 93  \\ 
140 & 107 & 5 & 37 & 31  \\ 
141 & 56 & 13 & 38 & 17  \\ 
142 & 229 & 25 & 197 & 55  \\ 
143 & 7 & 1 & 4 & 2  \\ 
144 & 41 & 12 & 23 & 14  \\ 
145 & 7 & 0 & 3 & 2  \\ 
146 & 49 & 16 & 9 & 18  \\ 
147 & 192 & 79 & 102 & 83  \\ 
148 & 59 & 1 & 53 & 13  \\ 
149 & 83 & 7 & 73 & 19  \\ 
150 & 7 & 0 & 1 & 2  \\ 
151 & 176 & 21 & 126 & 47  \\ 
152 & 67 & 7 & 29 & 19  \\ 
153 & 11 & 1 & 7 & 3  \\ 
154 & 419 & 94 & 197 & 128  \\ 
155 & 556 & 55 & 498 & 123  \\ 
156 & 64 & 1 & 44 & 17  \\ 
157 & 309 & 79 & 297 & 83  \\ 
158 & 331 & 34 & 173 & 92  \\ 
159 & 41 & 12 & 11 & 14  \\ 
160 & 7 & 1 & 1 & 2  \\ 
161 & 333 & 42 & 311 & 68  \\ 
162 & 417 & 34 & 183 & 116  \\ 
163 & 237 & 8 & 89 & 66  \\ 
164 & 166 & 11 & 146 & 37  \\ 
165 & 38 & 7 & 16 & 11  \\ 
166 & 243 & 37 & 77 & 69  \\ 
167 & 554 & 149 & 448 & 169  \\ 
168 & 34 & 1 & 22 & 9  \\ 
169 & 31 & 3 & 27 & 7  \\ 
170 & 29 & 7 & 3 & 9  \\ 
171 & 46 & 7 & 8 & 13  \\ 
172 & 982 & 33 & 222 & 271  \\ 
173 & 241 & 38 & 77 & 68  \\ 
174 & 13 & 1 & 11 & 3  \\ 
175 & 4 & 0 & 3 & 1  \\ 
176 & 21 & 6 & 1 & 7  \\ 
177 & 88 & 23 & 29 & 28  \\ 
178 & 67 & 11 & 53 & 17  \\ 
179 & 211 & 14 & 147 & 54  \\ 
180 & 38 & 7 & 34 & 9  \\ 
181 & 163 & 9 & 161 & 21  \\ 
182 & 103 & 26 & 47 & 32  \\ 
183 & 11 & 1 & 1 & 3  \\ 
184 & 475 & 135 & 123 & 157  \\ 
185 & 81 & 9 & 67 & 19  \\ 
186 & 43 & 12 & 19 & 14  \\ 
187 & 21 & 4 & 1 & 6  \\ 
188 & 218 & 43 & 158 & 59  \\ 
189 & 39 & 7 & 3 & 11  \\ 
190 & 79 & 10 & 73 & 16  \\ 
191 & 599 & 39 & 356 & 156  \\ 
192 & 62 & 9 & 2 & 17  \\ 
193 & 43 & 12 & 9 & 14  \\ 
194 & 233 & 46 & 209 & 56  \\ 
195 & 19 & 2 & 17 & 4  \\ 
196 & 57 & 16 & 41 & 18  \\ 
197 & 41 & 11 & 17 & 13  \\ 
198 & 71 & 7 & 17 & 19  \\ 
199 & 202 & 31 & 196 & 37  \\ 
200 & 266 & 13 & 262 & 35  \\ 
201 & 70 & 11 & 64 & 15  \\ 
202 & 757 & 299 & 253 & 313  \\ 
203 & 41 & 11 & 1 & 13  \\ 
204 & 59 & 2 & 53 & 12  \\ 
205 & 19 & 1 & 7 & 5  \\ 
206 & 4747 & 506 & 3923 & 1084  \\ 
207 & 36 & 2 & 33 & 7  \\ 
208 & 12 & 1 & 8 & 3  \\ 
209 & 71 & 14 & 5 & 20  \\ 
210 & 106 & 5 & 62 & 27  \\ 
211 & 1292 & 461 & 818 & 487  \\ 
212 & 226 & 77 & 14 & 83  \\ 
213 & 83 & 13 & 59 & 21  \\ 
214 & 2677 & 33 & 1821 & 659  \\ 
215 & 99 & 33 & 73 & 35  \\ 
216 & 42 & 11 & 18 & 13  \\ 
217 & 38 & 9 & 24 & 11  \\ 
218 & 127 & 29 & 49 & 37  \\ 
219 & 89 & 1 & 67 & 21  \\ 
220 & 114 & 5 & 106 & 21  \\ 
221 & 14 & 1 & 12 & 3  \\ 
222 & 31 & 7 & 1 & 9  \\ 
223 & 1399 & 353 & 1054 & 406  \\ 
224 & 414 & 3 & 62 & 107  \\ 
225 & 11 & 2 & 2 & 3  \\ 
226 & 723 & 4 & 229 & 186  \\ 
227 & 459 & 99 & 449 & 103  \\ 
228 & 212 & 27 & 136 & 53  \\ 
229 & 17 & 2 & 13 & 4  \\ 
230 & 556 & 25 & 464 & 121  \\ 
231 & 345 & 67 & 51 & 95  \\ 
232 & 139 & 23 & 43 & 37  \\ 
233 & 37 & 2 & 31 & 8  \\ 
234 & 111 & 23 & 33 & 31  \\ 
235 & 176 & 37 & 82 & 49  \\ 
236 & 1269 & 237 & 619 & 341  \\ 
237 & 82 & 5 & 76 & 15  \\ 
238 & 97 & 11 & 71 & 23  \\ 
239 & 120 & 0 & 119 & 13  \\ 
240 & 4 & 0 & 2 & 1  \\ 
241 & 56 & 13 & 31 & 16  \\ 
242 & 407 & 241 & 33 & 243  \\ 
243 & 12 & 1 & 6 & 3  \\ 
244 & 73 & 7 & 71 & 11  \\ 
245 & 99 & 27 & 1 & 31  \\ 
246 & 373 & 111 & 251 & 121  \\ 
247 & 782 & 91 & 244 & 199  \\ 
248 & 96 & 15 & 28 & 25  \\ 
249 & 689 & 37 & 473 & 163  \\ 
250 & 90 & 31 & 10 & 33  \\ 
251 & 334 & 73 & 168 & 93  \\ 
252 & 32 & 7 & 4 & 9  \\ 
253 & 4061 & 468 & 3969 & 614  \\ 
254 & 286 & 127 & 222 & 129  \\ 
255 & 23 & 4 & 11 & 6  \\ 
257 & 16 & 1 & 1 & 4  \\ 
258 & 549 & 1 & 483 & 109  \\ 
259 & 12 & 1 & 2 & 3  \\ 
260 & 9 & 0 & 7 & 2  \\ 
261 & 54 & 7 & 48 & 11  \\ 
262 & 3527 & 796 & 383 & 998  \\ 
263 & 289 & 121 & 237 & 123  \\ 
264 & 896 & 119 & 380 & 225  \\ 
265 & 33 & 7 & 17 & 9  \\ 
266 & 372 & 207 & 8 & 209  \\ 
267 & 297 & 41 & 237 & 67  \\ 
268 & 924 & 47 & 684 & 209  \\ 
269 & 23 & 4 & 3 & 6  \\ 
270 & 199 & 15 & 71 & 49  \\ 
271 & 303 & 104 & 239 & 108  \\ 
272 & 14 & 2 & 12 & 3  \\ 
273 & 8 & 1 & 1 & 2  \\ 
274 & 118 & 3 & 18 & 29  \\ 
275 & 13 & 1 & 9 & 3  \\ 
276 & 47 & 12 & 1 & 14  \\ 
277 & 87 & 1 & 43 & 21  \\ 
278 & 2119 & 592 & 1217 & 658  \\ 
279 & 252 & 49 & 66 & 67  \\ 
280 & 164 & 63 & 108 & 65  \\ 
281 & 188 & 23 & 152 & 41  \\ 
282 & 466 & 49 & 278 & 111  \\ 
283 & 357 & 184 & 209 & 186  \\ 
284 & 3306 & 687 & 3022 & 769  \\ 
285 & 47 & 1 & 29 & 11  \\ 
286 & 37 & 2 & 29 & 8  \\ 
287 & 17 & 2 & 11 & 4  \\ 
288 & 33 & 1 & 9 & 8  \\ 
289 & 31 & 2 & 22 & 7  \\ 
290 & 41 & 4 & 1 & 10  \\ 
291 & 364 & 19 & 248 & 83  \\ 
292 & 583 & 2 & 439 & 128  \\ 
293 & 117 & 42 & 31 & 44  \\ 
294 & 373 & 113 & 19 & 123  \\ 
295 & 191 & 29 & 163 & 41  \\ 
296 & 34 & 7 & 18 & 9  \\ 
297 & 54 & 1 & 12 & 13  \\ 
298 & 629 & 53 & 331 & 149  \\ 
299 & 438 & 33 & 436 & 43  \\ 
300 & 31 & 6 & 17 & 8  \\ 
301 & 64 & 13 & 22 & 17  \\ 
302 & 1233 & 123 & 881 & 277  \\ 
303 & 28 & 6 & 23 & 7  \\ 
304 & 196 & 13 & 32 & 47  \\ 
305 & 21 & 3 & 13 & 5  \\ 
306 & 104 & 9 & 76 & 23  \\ 
307 & 3792 & 569 & 722 & 939  \\ 
308 & 389 & 73 & 213 & 99  \\ 
309 & 277 & 8 & 71 & 66  \\ 
310 & 39 & 1 & 23 & 9  \\ 
311 & 433 & 88 & 189 & 114  \\ 
312 & 38 & 1 & 14 & 9  \\ 
313 & 51 & 2 & 23 & 12  \\ 
314 & 126 & 3 & 74 & 29  \\ 
315 & 79 & 15 & 61 & 19  \\ 
316 & 87 & 8 & 71 & 18  \\ 
317 & 203 & 46 & 1 & 56  \\ 
318 & 441 & 52 & 33 & 106  \\ 
319 & 12 & 2 & 1 & 3  \\ 
320 & 309 & 66 & 37 & 83  \\ 
321 & 187 & 12 & 134 & 41  \\ 
322 & 113 & 23 & 71 & 29  \\ 
323 & 101 & 7 & 89 & 19  \\ 
324 & 27 & 4 & 21 & 6  \\ 
325 & 775 & 53 & 415 & 179  \\ 
326 & 829 & 214 & 149 & 244  \\ 
327 & 13 & 1 & 7 & 3  \\ 
328 & 88 & 7 & 76 & 17  \\ 
329 & 243 & 15 & 149 & 55  \\ 
330 & 2263 & 104 & 1889 & 450  \\ 
331 & 1599 & 224 & 387 & 386  \\ 
332 & 2455 & 653 & 1625 & 721  \\ 
333 & 137 & 6 & 119 & 26  \\ 
334 & 2147 & 193 & 1861 & 413  \\ 
335 & 74 & 18 & 7 & 21  \\ 
336 & 41 & 4 & 29 & 9  \\ 
337 & 16 & 3 & 9 & 4  \\ 
338 & 17 & 2 & 7 & 4  \\ 
339 & 16 & 1 & 14 & 3  \\ 
340 & 130 & 17 & 10 & 31  \\ 
341 & 354 & 9 & 328 & 59  \\ 
342 & 799 & 57 & 569 & 173  \\ 
343 & 56 & 1 & 14 & 13  \\ 
344 & 338 & 19 & 178 & 77  \\ 
345 & 394 & 27 & 158 & 91  \\ 
346 & 262 & 77 & 146 & 83  \\ 
347 & 2428 & 797 & 1042 & 841  \\ 
348 & 17 & 2 & 1 & 4  \\ 
349 & 121 & 26 & 53 & 32  \\ 
350 & 31 & 1 & 17 & 7  \\ 
351 & 9 & 1 & 6 & 2  \\ 
352 & 158 & 59 & 62 & 61  \\ 
353 & 59 & 8 & 7 & 14  \\ 
354 & 148 & 1 & 88 & 33  \\ 
355 & 93 & 1 & 49 & 21  \\ 
356 & 209 & 46 & 31 & 56  \\ 
357 & 13 & 1 & 1 & 3  \\ 
358 & 5741 & 453 & 4283 & 1209  \\ 
359 & 3043 & 406 & 547 & 718  \\ 
360 & 22 & 3 & 14 & 5  \\ 
361 & 1179 & 314 & 987 & 334  \\ 
362 & 1464 & 171 & 304 & 341  \\ 
363 & 451 & 113 & 77 & 129  \\ 
364 & 128 & 13 & 68 & 29  \\ 
365 & 87 & 9 & 59 & 19  \\ 
366 & 363 & 7 & 333 & 61  \\ 
367 & 1054 & 332 & 781 & 347  \\ 
368 & 29 & 7 & 17 & 8  \\ 
369 & 5 & 0 & 4 & 1  \\ 
370 & 73 & 20 & 31 & 22  \\ 
371 & 9 & 0 & 5 & 2  \\ 
372 & 242 & 43 & 118 & 59  \\ 
373 & 61 & 6 & 3 & 14  \\ 
374 & 97 & 2 & 31 & 22  \\ 
375 & 22 & 1 & 4 & 5  \\ 
376 & 63 & 16 & 31 & 18  \\ 
377 & 31 & 4 & 27 & 6  \\ 
378 & 207 & 13 & 39 & 47  \\ 
379 & 4581 & 524 & 3823 & 906  \\ 
380 & 164 & 11 & 64 & 37  \\ 
381 & 223 & 72 & 31 & 76  \\ 
382 & 9529 & 677 & 3799 & 2147  \\ 
383 & 692 & 142 & 309 & 177  \\ 
384 & 164 & 3 & 44 & 37  \\ 
385 & 34 & 5 & 1 & 8  \\ 
386 & 197 & 23 & 131 & 43  \\ 
387 & 113 & 5 & 59 & 25  \\ 
388 & 36 & 7 & 16 & 9  \\ 
389 & 597 & 66 & 191 & 136  \\ 
390 & 29 & 5 & 11 & 7  \\ 
391 & 498 & 49 & 54 & 113  \\ 
392 & 102 & 9 & 94 & 17  \\ 
393 & 847 & 118 & 323 & 196  \\ 
394 & 101 & 19 & 35 & 25  \\ 
395 & 239 & 24 & 77 & 54  \\ 
396 & 58 & 3 & 14 & 13  \\ 
397 & 287 & 49 & 181 & 67  \\ 
398 & 939 & 246 & 653 & 268  \\ 
399 & 748 & 79 & 734 & 99  \\ 
400 & 149 & 34 & 43 & 40  \\ 
401 & 61 & 8 & 17 & 14  \\ 
402 & 439 & 16 & 97 & 98  \\ 
403 & 463 & 13 & 157 & 103  \\ 
404 & 166 & 21 & 158 & 27  \\ 
405 & 9 & 0 & 3 & 2  \\ 
406 & 103 & 11 & 47 & 23  \\ 
407 & 101 & 4 & 97 & 14  \\ 
408 & 184 & 17 & 124 & 39  \\ 
409 & 353 & 42 & 237 & 76  \\ 
410 & 9 & 0 & 1 & 2  \\ 
411 & 29 & 4 & 23 & 6  \\ 
412 & 623 & 189 & 417 & 199  \\ 
413 & 713 & 76 & 349 & 158  \\ 
414 & 462 & 37 & 366 & 91  \\ 
415 & 97 & 6 & 69 & 20  \\ 
416 & 77 & 1 & 53 & 16  \\ 
417 & 301 & 134 & 23 & 136  \\ 
418 & 1394 & 179 & 938 & 301  \\ 
419 & 1689 & 114 & 851 & 368  \\ 
420 & 79 & 7 & 47 & 17  \\ 
421 & 996 & 31 & 586 & 213  \\ 
422 & 239 & 22 & 183 & 48  \\ 
423 & 211 & 21 & 23 & 47  \\ 
424 & 19 & 2 & 13 & 4  \\ 
425 & 23 & 1 & 11 & 5  \\ 
426 & 439 & 23 & 271 & 93  \\ 
427 & 52 & 7 & 38 & 11  \\ 
428 & 291 & 18 & 77 & 64  \\ 
429 & 14 & 1 & 8 & 3  \\ 
430 & 63 & 11 & 23 & 15  \\ 
431 & 802 & 112 & 371 & 181  \\ 
432 & 9 & 1 & 3 & 2  \\ 
433 & 34 & 3 & 24 & 7  \\ 
434 & 174 & 63 & 50 & 65  \\ 
435 & 101 & 30 & 23 & 32  \\ 
436 & 174 & 23 & 74 & 39  \\ 
437 & 3831 & 907 & 1551 & 1037  \\ 
438 & 142 & 1 & 38 & 31  \\ 
439 & 541 & 55 & 337 & 115  \\ 
440 & 213 & 80 & 117 & 82  \\ 
441 & 482 & 27 & 302 & 101  \\ 
442 & 59 & 14 & 19 & 16  \\ 
443 & 4093 & 8 & 2109 & 876  \\ 
444 & 149 & 17 & 43 & 33  \\ 
445 & 376 & 3 & 324 & 67  \\ 
446 & 2796 & 393 & 1012 & 631  \\ 
447 & 37 & 7 & 11 & 9  \\ 
448 & 219 & 32 & 117 & 49  \\ 
449 & 153 & 36 & 76 & 41  \\ 
450 & 207 & 8 & 201 & 26  \\ 
451 & 69 & 16 & 47 & 18  \\ 
452 & 129 & 6 & 127 & 14  \\ 
453 & 169 & 25 & 133 & 35  \\ 
454 & 331 & 147 & 123 & 149  \\ 
455 & 33 & 3 & 19 & 7  \\ 
456 & 862 & 99 & 658 & 173  \\ 
457 & 259 & 22 & 103 & 56  \\ 
458 & 197 & 19 & 11 & 43  \\ 
459 & 63 & 8 & 3 & 14  \\ 
460 & 277 & 19 & 139 & 59  \\ 
461 & 791 & 157 & 627 & 183  \\ 
462 & 326 & 39 & 158 & 71  \\ 
463 & 277 & 33 & 186 & 58  \\ 
464 & 14 & 1 & 6 & 3  \\ 
465 & 97 & 8 & 4 & 21  \\ 
466 & 127 & 7 & 119 & 19  \\ 
467 & 576 & 3 & 358 & 119  \\ 
468 & 388 & 1 & 224 & 81  \\ 
469 & 461 & 8 & 209 & 98  \\ 
470 & 1139 & 109 & 11 & 247  \\ 
471 & 127 & 25 & 55 & 31  \\ 
472 & 484 & 67 & 224 & 107  \\ 
473 & 3838 & 109 & 2784 & 759  \\ 
474 & 1129 & 38 & 23 & 242  \\ 
475 & 61 & 1 & 23 & 13  \\ 
476 & 38 & 7 & 18 & 9  \\ 
477 & 14 & 1 & 4 & 3  \\ 
478 & 267 & 23 & 211 & 51  \\ 
479 & 645 & 30 & 313 & 136  \\ 
480 & 14 & 1 & 2 & 3  \\ 
481 & 33 & 8 & 4 & 9  \\ 
482 & 19 & 2 & 11 & 4  \\ 
483 & 1933 & 32 & 461 & 412  \\ 
484 & 362 & 19 & 122 & 77  \\ 
485 & 106 & 19 & 4 & 25  \\ 
486 & 39 & 3 & 33 & 7  \\ 
487 & 12817 & 2097 & 5689 & 2919  \\ 
488 & 174 & 3 & 38 & 37  \\ 
489 & 232 & 3 & 94 & 49  \\ 
490 & 1421 & 254 & 371 & 334  \\ 
491 & 3982 & 1513 & 2892 & 1539  \\ 
492 & 38 & 7 & 14 & 9  \\ 
493 & 18 & 1 & 16 & 3  \\ 
494 & 108 & 9 & 32 & 23  \\ 
495 & 19 & 1 & 8 & 4  \\ 
496 & 36 & 8 & 26 & 9  \\ 
497 & 109 & 2 & 38 & 23  \\ 
498 & 383 & 21 & 217 & 79  \\ 
499 & 2361 & 577 & 365 & 645  \\ 
500 & 509 & 80 & 237 & 114  \\ 
501 & 4666 & 949 & 1012 & 1151  \\ 
502 & 479 & 68 & 23 & 106  \\ 
503 & 792 & 363 & 214 & 367  \\ 
504 & 104 & 13 & 76 & 21  \\ 
505 & 11 & 0 & 9 & 2  \\ 
506 & 67 & 16 & 21 & 18  \\ 
507 & 428 & 41 & 122 & 91  \\ 
508 & 802 & 99 & 294 & 173  \\ 
509 & 105 & 15 & 47 & 23  \\ 
510 & 23 & 3 & 7 & 5  \\ 
511 & 19 & 1 & 2 & 4  \\ 
512 & 278 & 17 & 122 & 58  \\ 
513 & 66 & 13 & 9 & 16  \\ 
514 & 541 & 21 & 339 & 109  \\ 
515 & 1827 & 24 & 797 & 380  \\ 
516 & 176 & 11 & 4 & 37  \\ 
517 & 257 & 12 & 213 & 46  \\ 
518 & 162 & 3 & 134 & 29  \\ 
519 & 59 & 7 & 49 & 11  \\ 
520 & 28 & 3 & 24 & 5  \\ 
521 & 53 & 2 & 49 & 8  \\ 
522 & 281 & 18 & 223 & 52  \\ 
523 & 1253 & 93 & 207 & 263  \\ 
524 & 491 & 71 & 229 & 107  \\ 
525 & 3343 & 371 & 907 & 711  \\ 
526 & 5659 & 372 & 1179 & 1184  \\ 
527 & 24 & 0 & 7 & 5  \\ 
528 & 14 & 2 & 8 & 3  \\ 
529 & 2074 & 157 & 487 & 434  \\ 
530 & 149 & 1 & 43 & 31  \\ 
531 & 216 & 7 & 138 & 43  \\ 
532 & 639 & 29 & 159 & 133  \\ 
533 & 19 & 2 & 7 & 4  \\ 
534 & 183 & 16 & 81 & 38  \\ 
535 & 9836 & 1161 & 5786 & 2037  \\ 
536 & 2307 & 372 & 29 & 518  \\ 
537 & 543 & 7 & 531 & 61  \\ 
538 & 718 & 137 & 358 & 169  \\ 
539 & 743 & 77 & 335 & 155  \\ 
540 & 27 & 1 & 21 & 5  \\ 
541 & 104 & 29 & 2 & 31  \\ 
542 & 3259 & 622 & 1091 & 772  \\ 
543 & 19 & 2 & 1 & 4  \\ 
544 & 5 & 0 & 3 & 1  \\ 
545 & 42 & 3 & 31 & 8  \\ 
546 & 236 & 17 & 16 & 49  \\ 
547 & 921 & 82 & 173 & 192  \\ 
548 & 252 & 27 & 8 & 53  \\ 
549 & 135 & 10 & 21 & 28  \\ 
550 & 233 & 46 & 31 & 56  \\ 
551 & 1122 & 27 & 514 & 229  \\ 
552 & 776 & 23 & 316 & 159  \\ 
553 & 853 & 65 & 379 & 175  \\ 
554 & 1105 & 180 & 457 & 246  \\ 
555 & 83 & 21 & 17 & 23  \\ 
556 & 2310 & 11 & 642 & 475  \\ 
557 & 127 & 8 & 53 & 26  \\ 
558 & 763 & 21 & 43 & 157  \\ 
559 & 404 & 8 & 103 & 83  \\ 
560 & 16 & 1 & 12 & 3  \\ 
561 & 128 & 3 & 59 & 26  \\ 
562 & 229 & 16 & 227 & 22  \\ 
563 & 1739 & 256 & 513 & 378  \\ 
564 & 683 & 22 & 539 & 124  \\ 
565 & 347 & 19 & 121 & 71  \\ 
566 & 618 & 51 & 514 & 109  \\ 
567 & 39 & 2 & 3 & 8  \\ 
568 & 499 & 71 & 211 & 107  \\ 
569 & 92 & 17 & 54 & 21  \\ 
570 & 1994 & 163 & 1958 & 227  \\ 
571 & 2033 & 18 & 891 & 412  \\ 
572 & 367 & 5 & 81 & 75  \\ 
573 & 2269 & 547 & 359 & 607  \\ 
574 & 69 & 16 & 13 & 18  \\ 
575 & 49 & 2 & 43 & 8  \\ 
576 & 34 & 3 & 2 & 7  \\ 
577 & 164 & 43 & 103 & 46  \\ 
578 & 152 & 1 & 4 & 31  \\ 
579 & 67 & 8 & 23 & 14  \\ 
580 & 1360 & 371 & 200 & 397  \\ 
581 & 709 & 22 & 453 & 138  \\ 
582 & 301 & 37 & 107 & 63  \\ 
583 & 139 & 26 & 73 & 32  \\ 
584 & 54 & 3 & 2 & 11  \\ 
585 & 19 & 1 & 17 & 3  \\ 
586 & 109 & 9 & 99 & 17  \\ 
587 & 1867 & 71 & 481 & 379  \\ 
588 & 214 & 11 & 178 & 37  \\ 
589 & 112 & 29 & 74 & 31  \\ 
590 & 292 & 77 & 56 & 83  \\ 
591 & 14 & 2 & 1 & 3  \\ 
592 & 52 & 7 & 16 & 11  \\ 
593 & 107 & 1 & 63 & 21  \\ 
594 & 69 & 4 & 3 & 14  \\ 
595 & 178 & 29 & 172 & 31  \\ 
596 & 96 & 3 & 52 & 19  \\ 
597 & 1694 & 449 & 296 & 483  \\ 
598 & 113 & 8 & 71 & 22  \\ 
599 & 1296 & 93 & 98 & 263  \\ 
600 & 468 & 41 & 168 & 95  \\ 
601 & 1124 & 7 & 134 & 227  \\ 
602 & 127 & 23 & 41 & 29  \\ 
603 & 496 & 37 & 442 & 79  \\ 
604 & 2824 & 937 & 1012 & 967  \\ 
605 & 394 & 5 & 152 & 79  \\ 
606 & 199 & 23 & 71 & 41  \\ 
607 & 653 & 421 & 46 & 422  \\ 
608 & 144 & 3 & 8 & 29  \\ 
609 & 5 & 0 & 2 & 1  \\ 
610 & 241 & 32 & 119 & 50  \\ 
611 & 61 & 2 & 33 & 12  \\ 
612 & 41 & 7 & 23 & 9  \\ 
613 & 308 & 1 & 254 & 53  \\ 
614 & 4261 & 68 & 37 & 856  \\ 
615 & 94 & 17 & 56 & 21  \\ 
616 & 666 & 29 & 258 & 133  \\ 
617 & 5366 & 1003 & 1726 & 1237  \\ 
618 & 501 & 98 & 117 & 118  \\ 
619 & 1242 & 11 & 4 & 249  \\ 
620 & 539 & 80 & 267 & 114  \\ 
621 & 657 & 73 & 381 & 131  \\ 
622 & 6937 & 742 & 5071 & 1312  \\ 
623 & 23 & 2 & 19 & 4  \\ 
624 & 41 & 2 & 23 & 8  \\ 
626 & 25 & 1 & 1 & 5  \\ 
627 & 101 & 13 & 13 & 21  \\ 
628 & 223 & 59 & 193 & 61  \\ 
629 & 101 & 26 & 47 & 28  \\ 
630 & 309 & 7 & 141 & 61  \\ 
631 & 796 & 283 & 466 & 289  \\ 
632 & 464 & 91 & 148 & 109  \\ 
633 & 1277 & 64 & 433 & 254  \\ 
634 & 53 & 7 & 19 & 11  \\ 
635 & 388 & 5 & 374 & 47  \\ 
636 & 124 & 7 & 88 & 23  \\ 
637 & 244 & 13 & 146 & 47  \\ 
638 & 404 & 19 & 256 & 77  \\ 
639 & 314 & 7 & 172 & 61  \\ 
640 & 29 & 6 & 3 & 7  \\ 
641 & 25 & 2 & 4 & 5  \\ 
642 & 467 & 27 & 389 & 79  \\ 
643 & 10264 & 557 & 1262 & 2041  \\ 
644 & 1294 & 111 & 650 & 255  \\ 
645 & 113 & 4 & 59 & 22  \\ 
646 & 1103 & 148 & 569 & 226  \\ 
647 & 688 & 203 & 606 & 207  \\ 
648 & 696 & 61 & 204 & 139  \\ 
649 & 217 & 7 & 19 & 43  \\ 
650 & 567 & 49 & 183 & 113  \\ 
651 & 197 & 32 & 113 & 42  \\ 
652 & 4073 & 127 & 3751 & 587  \\ 
653 & 629 & 131 & 311 & 151  \\ 
654 & 1073 & 31 & 529 & 209  \\ 
655 & 2509 & 105 & 504 & 496  \\ 
656 & 96 & 9 & 44 & 19  \\ 
657 & 14 & 1 & 13 & 2  \\ 
658 & 3811 & 187 & 2299 & 727  \\ 
659 & 6067 & 1342 & 1841 & 1516  \\ 
660 & 439 & 93 & 23 & 107  \\ 
661 & 789 & 104 & 313 & 162  \\ 
662 & 4291 & 811 & 1643 & 983  \\ 
663 & 41 & 4 & 23 & 8  \\ 
664 & 4360 & 125 & 44 & 859  \\ 
665 & 1097 & 92 & 1069 & 130  \\ 
666 & 29 & 4 & 11 & 6  \\ 
667 & 37 & 6 & 9 & 8  \\ 
668 & 496 & 51 & 172 & 99  \\ 
669 & 529 & 2 & 83 & 104  \\ 
670 & 1063 & 39 & 9 & 209  \\ 
671 & 6 & 0 & 5 & 1  \\ 
672 & 83 & 1 & 43 & 16  \\ 
673 & 47 & 6 & 34 & 9  \\ 
674 & 107 & 1 & 3 & 21  \\ 
675 & 75 & 17 & 15 & 19  \\ 
676 & 368 & 13 & 308 & 61  \\ 
677 & 811 & 23 & 99 & 159  \\ 
678 & 59 & 12 & 19 & 14  \\ 
679 & 479 & 14 & 353 & 86  \\ 
680 & 29 & 4 & 3 & 6  \\ 
681 & 1028 & 511 & 788 & 513  \\ 
682 & 399 & 28 & 151 & 78  \\ 
683 & 752 & 11 & 614 & 127  \\ 
684 & 688 & 71 & 148 & 137  \\ 
685 & 29 & 1 & 23 & 5  \\ 
686 & 63 & 3 & 49 & 11  \\ 
687 & 17 & 1 & 13 & 3  \\ 
688 & 734 & 257 & 46 & 263  \\ 
689 & 41 & 2 & 11 & 8  \\ 
690 & 374 & 43 & 98 & 75  \\ 
691 & 4309 & 191 & 163 & 841  \\ 
692 & 41 & 2 & 7 & 8  \\ 
693 & 54 & 7 & 12 & 11  \\ 
694 & 22222 & 4069 & 8342 & 4987  \\ 
695 & 241 & 13 & 37 & 47  \\ 
696 & 83 & 12 & 77 & 14  \\ 
697 & 192 & 19 & 158 & 33  \\ 
698 & 49 & 1 & 33 & 9  \\ 
699 & 427 & 89 & 257 & 101  \\ 
700 & 366 & 17 & 362 & 33  \\ 
701 & 31 & 4 & 21 & 6  \\ 
702 & 106 & 11 & 2 & 21  \\ 
703 & 11 & 1 & 8 & 2  \\ 
704 & 92 & 13 & 4 & 19  \\ 
705 & 211 & 11 & 23 & 41  \\ 
706 & 25 & 3 & 9 & 5  \\ 
707 & 242 & 13 & 4 & 47  \\ 
708 & 902 & 177 & 278 & 209  \\ 
709 & 97 & 2 & 95 & 10  \\ 
710 & 277 & 76 & 149 & 80  \\ 
711 & 291 & 79 & 183 & 83  \\ 
712 & 294 & 17 & 62 & 57  \\ 
713 & 784 & 43 & 412 & 149  \\ 
714 & 117 & 26 & 93 & 28  \\ 
715 & 53 & 4 & 31 & 10  \\ 
716 & 1847 & 401 & 57 & 453  \\ 
717 & 491 & 25 & 13 & 95  \\ 
718 & 1694 & 213 & 258 & 341  \\ 
719 & 1354 & 170 & 803 & 265  \\ 
720 & 19 & 2 & 17 & 3  \\ 
721 & 124 & 11 & 82 & 23  \\ 
722 & 493 & 16 & 229 & 94  \\ 
723 & 149 & 1 & 137 & 21  \\ 
724 & 37 & 3 & 21 & 7  \\ 
725 & 57 & 3 & 1 & 11  \\ 
726 & 1993 & 374 & 1153 & 444  \\ 
727 & 8989 & 1448 & 9913 & 554  \\ 
728 & 79 & 17 & 47 & 19  \\ 
729 & 87 & 9 & 33 & 17  \\ 
730 & 246 & 25 & 218 & 39  \\ 
731 & 58 & 1 & 28 & 11  \\ 
732 & 73 & 2 & 71 & 8  \\ 
733 & 7207 & 942 & 2027 & 1452  \\ 
734 & 1971 & 591 & 1237 & 611  \\ 
735 & 67 & 7 & 31 & 13  \\ 
736 & 234 & 59 & 134 & 63  \\ 
737 & 193 & 40 & 126 & 45  \\ 
738 & 49 & 1 & 31 & 9  \\ 
739 & 1682 & 101 & 1274 & 293  \\ 
740 & 92 & 1 & 56 & 17  \\ 
741 & 137 & 4 & 61 & 26  \\ 
742 & 73 & 4 & 17 & 14  \\ 
743 & 15959 & 2254 & 9303 & 3186  \\ 
744 & 446 & 1 & 422 & 57  \\ 
745 & 224 & 35 & 22 & 47  \\ 
746 & 1619 & 149 & 1339 & 271  \\ 
747 & 327 & 28 & 171 & 62  \\ 
748 & 42 & 7 & 2 & 9  \\ 
749 & 716 & 73 & 502 & 131  \\ 
750 & 194 & 5 & 58 & 37  \\ 
751 & 2399 & 1119 & 1356 & 1126  \\ 
752 & 48 & 6 & 46 & 7  \\ 
753 & 157 & 8 & 94 & 29  \\ 
754 & 136 & 9 & 108 & 23  \\ 
755 & 319 & 19 & 17 & 61  \\ 
756 & 177 & 8 & 159 & 26  \\ 
757 & 1578 & 133 & 488 & 303  \\ 
758 & 3101 & 43 & 69 & 591  \\ 
759 & 182 & 27 & 94 & 37  \\ 
760 & 606 & 1 & 454 & 105  \\ 
761 & 40 & 1 & 38 & 5  \\ 
762 & 4145 & 444 & 335 & 808  \\ 
763 & 103 & 4 & 89 & 16  \\ 
764 & 839 & 288 & 457 & 294  \\ 
765 & 11 & 0 & 7 & 2  \\ 
766 & 1161 & 44 & 1137 & 118  \\ 
767 & 38 & 8 & 21 & 9  \\ 
768 & 16 & 1 & 8 & 3  \\ 
769 & 37 & 1 & 13 & 7  \\ 
770 & 89 & 7 & 23 & 17  \\ 
771 & 724 & 81 & 238 & 141  \\ 
772 & 113 & 3 & 111 & 11  \\ 
773 & 2099 & 56 & 943 & 394  \\ 
774 & 507 & 131 & 267 & 139  \\ 
775 & 596 & 75 & 272 & 117  \\ 
776 & 358 & 35 & 30 & 69  \\ 
777 & 361 & 3 & 353 & 37  \\ 
778 & 754 & 21 & 426 & 139  \\ 
779 & 79 & 5 & 3 & 15  \\ 
780 & 160 & 17 & 20 & 31  \\ 
781 & 373 & 29 & 89 & 71  \\ 
782 & 877 & 136 & 181 & 182  \\ 
783 & 16 & 2 & 11 & 3  \\ 
784 & 88 & 17 & 10 & 20  \\ 
785 & 63 & 12 & 37 & 14  \\ 
786 & 1353 & 244 & 567 & 296  \\ 
787 & 4388 & 93 & 2814 & 791  \\ 
788 & 191 & 56 & 113 & 58  \\ 
789 & 891 & 43 & 687 & 151  \\ 
790 & 314 & 77 & 2 & 83  \\ 
791 & 843 & 96 & 11 & 164  \\ 
792 & 646 & 67 & 542 & 107  \\ 
793 & 201 & 16 & 137 & 36  \\ 
794 & 426 & 51 & 158 & 83  \\ 
795 & 109 & 11 & 103 & 15  \\ 
796 & 231 & 35 & 167 & 45  \\ 
797 & 513 & 64 & 309 & 98  \\ 
798 & 761 & 89 & 607 & 133  \\ 
799 & 24 & 2 & 23 & 3  \\ 
800 & 209 & 2 & 87 & 39  \\ 
801 & 31 & 6 & 14 & 7  \\ 
802 & 21 & 1 & 19 & 3  \\ 
803 & 16 & 1 & 6 & 3  \\ 
804 & 877 & 158 & 61 & 192  \\ 
805 & 37 & 3 & 9 & 7  \\ 
806 & 43 & 2 & 19 & 8  \\ 
807 & 1863 & 248 & 1371 & 346  \\ 
808 & 386 & 31 & 22 & 73  \\ 
809 & 5106 & 453 & 1870 & 965  \\ 
810 & 21 & 2 & 3 & 4  \\ 
811 & 4638 & 791 & 3016 & 963  \\ 
812 & 72 & 3 & 44 & 13  \\ 
813 & 1668 & 157 & 1584 & 221  \\ 
814 & 927 & 89 & 553 & 171  \\ 
815 & 127 & 11 & 36 & 24  \\ 
816 & 11 & 1 & 7 & 2  \\ 
817 & 31 & 6 & 12 & 7  \\ 
818 & 23 & 2 & 17 & 4  \\ 
819 & 16 & 1 & 2 & 3  \\ 
820 & 49 & 8 & 31 & 10  \\ 
821 & 39 & 6 & 11 & 8  \\ 
822 & 491 & 71 & 83 & 99  \\ 
823 & 982 & 259 & 664 & 271  \\ 
824 & 2314 & 289 & 982 & 449  \\ 
825 & 91 & 15 & 41 & 19  \\ 
826 & 1467 & 67 & 1011 & 257  \\ 
827 & 891 & 9 & 763 & 137  \\ 
828 & 499 & 23 & 131 & 93  \\ 
829 & 145 & 25 & 23 & 31  \\ 
830 & 1337 & 23 & 673 & 245  \\ 
831 & 32 & 6 & 19 & 7  \\ 
832 & 18 & 1 & 14 & 3  \\ 
833 & 122 & 7 & 73 & 22  \\ 
834 & 2083 & 3 & 1805 & 315  \\ 
835 & 8165 & 1301 & 4825 & 1657  \\ 
836 & 651 & 68 & 603 & 94  \\ 
837 & 168 & 37 & 18 & 41  \\ 
838 & 3773 & 416 & 421 & 722  \\ 
839 & 1121 & 76 & 557 & 206  \\ 
840 & 26 & 3 & 2 & 5  \\ 
841 & 1839 & 208 & 1391 & 324  \\ 
842 & 1259 & 131 & 1147 & 187  \\ 
843 & 854 & 21 & 832 & 89  \\ 
844 & 2303 & 479 & 1073 & 539  \\ 
845 & 91 & 7 & 13 & 17  \\ 
846 & 261 & 23 & 21 & 49  \\ 
847 & 134 & 10 & 13 & 25  \\ 
848 & 31 & 1 & 29 & 4  \\ 
849 & 603 & 172 & 246 & 179  \\ 
850 & 97 & 20 & 71 & 22  \\ 
851 & 63 & 12 & 29 & 14  \\ 
852 & 2098 & 811 & 962 & 821  \\ 
853 & 503 & 104 & 311 & 116  \\ 
854 & 31 & 1 & 25 & 5  \\ 
855 & 216 & 23 & 12 & 41  \\ 
856 & 1286 & 53 & 1282 & 83  \\ 
857 & 108 & 9 & 74 & 19  \\ 
858 & 172 & 3 & 128 & 29  \\ 
859 & 1097 & 141 & 621 & 209  \\ 
860 & 979 & 209 & 397 & 233  \\ 
861 & 17 & 1 & 11 & 3  \\ 
862 & 8767 & 712 & 4457 & 1606  \\ 
863 & 601 & 331 & 262 & 332  \\ 
864 & 117 & 2 & 93 & 19  \\ 
865 & 71 & 1 & 29 & 13  \\ 
866 & 237 & 18 & 11 & 44  \\ 
867 & 173 & 49 & 73 & 51  \\ 
868 & 337 & 82 & 97 & 88  \\ 
869 & 2834 & 487 & 326 & 601  \\ 
870 & 13 & 0 & 11 & 2  \\ 
871 & 383 & 29 & 19 & 71  \\ 
872 & 226 & 71 & 46 & 73  \\ 
873 & 151 & 26 & 43 & 32  \\ 
874 & 99 & 8 & 53 & 18  \\ 
875 & 80 & 17 & 10 & 19  \\ 
876 & 11 & 0 & 5 & 2  \\ 
877 & 712 & 43 & 518 & 121  \\ 
878 & 631 & 82 & 247 & 122  \\ 
879 & 79 & 1 & 61 & 13  \\ 
880 & 28 & 1 & 16 & 5  \\ 
881 & 25 & 4 & 16 & 5  \\ 
882 & 203 & 47 & 77 & 51  \\ 
883 & 3049 & 166 & 1283 & 556  \\ 
884 & 83 & 11 & 77 & 13  \\ 
885 & 139 & 10 & 97 & 24  \\ 
886 & 1657 & 413 & 1001 & 437  \\ 
887 & 4019 & 592 & 1303 & 802  \\ 
888 & 31 & 4 & 1 & 6  \\ 
889 & 401 & 32 & 107 & 74  \\ 
890 & 443 & 67 & 157 & 89  \\ 
891 & 63 & 12 & 3 & 14  \\ 
892 & 2252 & 147 & 2208 & 227  \\ 
893 & 408 & 9 & 314 & 67  \\ 
894 & 79 & 6 & 49 & 14  \\ 
895 & 501 & 47 & 143 & 93  \\ 
896 & 44 & 7 & 12 & 9  \\ 
897 & 38 & 1 & 31 & 6  \\ 
898 & 477 & 1 & 421 & 69  \\ 
899 & 16388 & 203 & 6242 & 2977  \\ 
900 & 44 & 7 & 8 & 9  \\ 
901 & 31 & 2 & 29 & 4  \\ 
902 & 983 & 49 & 337 & 179  \\ 
903 & 601 & 72 & 517 & 98  \\ 
904 & 94 & 9 & 54 & 17  \\ 
905 & 79 & 8 & 53 & 14  \\ 
906 & 721 & 22 & 487 & 124  \\ 
907 & 9553 & 2087 & 4111 & 2297  \\ 
908 & 6359 & 382 & 2727 & 1152  \\ 
909 & 89 & 11 & 19 & 17  \\ 
910 & 11 & 0 & 3 & 2  \\ 
911 & 1199 & 226 & 623 & 262  \\ 
912 & 49 & 4 & 11 & 9  \\ 
913 & 659 & 72 & 254 & 123  \\ 
914 & 266 & 1 & 254 & 31  \\ 
915 & 11 & 0 & 1 & 2  \\ 
916 & 468 & 1 & 304 & 81  \\ 
917 & 404 & 103 & 142 & 109  \\ 
918 & 63 & 1 & 39 & 11  \\ 
919 & 1489 & 308 & 349 & 346  \\ 
920 & 1981 & 4 & 733 & 358  \\ 
921 & 1121 & 273 & 721 & 289  \\ 
922 & 613 & 75 & 309 & 115  \\ 
923 & 463 & 43 & 389 & 73  \\ 
924 & 172 & 19 & 136 & 29  \\ 
925 & 1321 & 95 & 103 & 241  \\ 
926 & 2383 & 514 & 1457 & 562  \\ 
927 & 171 & 17 & 138 & 28  \\ 
928 & 13 & 1 & 11 & 2  \\ 
929 & 89 & 2 & 37 & 16  \\ 
930 & 269 & 15 & 227 & 41  \\ 
931 & 637 & 11 & 427 & 109  \\ 
932 & 94 & 1 & 22 & 17  \\ 
933 & 563 & 47 & 59 & 103  \\ 
934 & 4895 & 1278 & 225 & 1346  \\ 
935 & 111 & 6 & 43 & 20  \\ 
936 & 64 & 11 & 44 & 13  \\ 
937 & 438 & 77 & 356 & 87  \\ 
938 & 709 & 2 & 173 & 128  \\ 
939 & 128 & 19 & 70 & 25  \\ 
940 & 1407 & 212 & 849 & 274  \\ 
941 & 39 & 3 & 19 & 7  \\ 
942 & 61 & 3 & 19 & 11  \\ 
943 & 32 & 6 & 9 & 7  \\ 
944 & 1308 & 93 & 1288 & 127  \\ 
945 & 12 & 1 & 9 & 2  \\ 
946 & 4991 & 326 & 2153 & 896  \\ 
947 & 1171 & 148 & 723 & 216  \\ 
948 & 1312 & 281 & 268 & 311  \\ 
949 & 866 & 63 & 168 & 157  \\ 
950 & 27 & 3 & 11 & 5  \\ 
951 & 47 & 7 & 31 & 9  \\ 
952 & 688 & 7 & 376 & 121  \\ 
953 & 127 & 1 & 93 & 21  \\ 
954 & 115 & 20 & 101 & 22  \\ 
955 & 328 & 9 & 54 & 59  \\ 
956 & 1978 & 101 & 1846 & 251  \\ 
957 & 254 & 27 & 32 & 47  \\ 
958 & 491 & 41 & 467 & 61  \\ 
959 & 11 & 1 & 4 & 2  \\ 
960 & 379 & 35 & 133 & 69  \\ 
961 & 1473 & 32 & 449 & 264  \\ 
962 & 17 & 1 & 9 & 3  \\ 
963 & 1137 & 259 & 147 & 281  \\ 
964 & 369 & 104 & 113 & 108  \\ 
965 & 211 & 35 & 107 & 43  \\ 
966 & 239 & 15 & 55 & 43  \\ 
967 & 252386 & 98431 & 416776 & 90427  \\ 
968 & 486 & 3 & 482 & 37  \\ 
969 & 28 & 1 & 10 & 5  \\ 
970 & 193 & 34 & 119 & 40  \\ 
971 & 1081 & 109 & 861 & 177  \\ 
972 & 93 & 9 & 3 & 17  \\ 
973 & 2186 & 389 & 796 & 463  \\ 
974 & 1511 & 544 & 537 & 552  \\ 
975 & 70 & 9 & 40 & 13  \\ 
976 & 11 & 1 & 1 & 2  \\ 
977 & 224 & 11 & 218 & 23  \\ 
978 & 703 & 103 & 601 & 123  \\ 
979 & 122 & 29 & 56 & 31  \\ 
980 & 757 & 50 & 517 & 128  \\ 
981 & 113 & 13 & 23 & 21  \\ 
982 & 3597 & 596 & 331 & 738  \\ 
983 & 5980 & 167 & 3850 & 1019  \\ 
984 & 242 & 19 & 118 & 43  \\ 
985 & 28 & 1 & 2 & 5  \\ 
986 & 417 & 6 & 161 & 74  \\ 
987 & 157 & 2 & 31 & 28  \\ 
988 & 4931 & 366 & 579 & 886  \\ 
989 & 1327 & 62 & 467 & 236  \\ 
990 & 161 & 20 & 73 & 30  \\ 
991 & 642 & 196 & 349 & 201  \\ 
992 & 466 & 1 & 94 & 83  \\ 
993 & 497 & 14 & 496 & 27  \\ 
994 & 4822 & 223 & 2834 & 833  \\ 
995 & 517 & 14 & 119 & 92  \\ 
996 & 3251 & 194 & 2227 & 546  \\ 
997 & 137 & 9 & 123 & 19  \\ 
998 & 687 & 27 & 311 & 121  \\ 
999 & 41 & 6 & 13 & 8  \\ 
\end{longtable}

\label{longtable}
\section*{Conclusions}
We can state a conjecture as follows.\\  
\textbf{Conjecture:}\\
Let $h$ be an arbitrary integral number.\\  
Then the equation $A^4 + hB^4 = C^4 + hD^4$ always has the integral solutions.
\vskip\baselineskip

\section*{Acknowledgement}
I would like to thank Jaroslaw Wroblewski for providing many solutions.
\vskip\baselineskip

\section*{Appendix} 
We show the link for the numerical solutions for $1000<h<20000$\\  
\url{http://www.maroon.dti.ne.jp/fermat/a^4+hb^4=c^4+d^4/1000.pdf}\\  
\url{http://www.maroon.dti.ne.jp/fermat/a^4+hb^4=c^4+d^4/2000.pdf}\\
\url{http://www.maroon.dti.ne.jp/fermat/a^4+hb^4=c^4+d^4/3000.pdf}\\
\url{http://www.maroon.dti.ne.jp/fermat/a^4+hb^4=c^4+d^4/4000.pdf}\\
\url{http://www.maroon.dti.ne.jp/fermat/a^4+hb^4=c^4+d^4/5000.pdf}\\
\url{http://www.maroon.dti.ne.jp/fermat/a^4+hb^4=c^4+d^4/6000.pdf}\\
\url{http://www.maroon.dti.ne.jp/fermat/a^4+hb^4=c^4+d^4/7000.pdf}\\
\url{http://www.maroon.dti.ne.jp/fermat/a^4+hb^4=c^4+d^4/8000.pdf}\\
\url{http://www.maroon.dti.ne.jp/fermat/a^4+hb^4=c^4+d^4/9000.pdf}\\
\url{http://www.maroon.dti.ne.jp/fermat/a^4+hb^4=c^4+d^4/10000.pdf}\\
\url{http://www.maroon.dti.ne.jp/fermat/a^4+hb^4=c^4+d^4/11000.pdf}\\
\url{http://www.maroon.dti.ne.jp/fermat/a^4+hb^4=c^4+d^4/12000.pdf}\\
\url{http://www.maroon.dti.ne.jp/fermat/a^4+hb^4=c^4+d^4/13000.pdf}\\
\url{http://www.maroon.dti.ne.jp/fermat/a^4+hb^4=c^4+d^4/14000.pdf}\\
\url{http://www.maroon.dti.ne.jp/fermat/a^4+hb^4=c^4+d^4/15000.pdf}\\
\url{http://www.maroon.dti.ne.jp/fermat/a^4+hb^4=c^4+d^4/16000.pdf}\\
\url{http://www.maroon.dti.ne.jp/fermat/a^4+hb^4=c^4+d^4/17000.pdf}\\
\url{http://www.maroon.dti.ne.jp/fermat/a^4+hb^4=c^4+d^4/18000.pdf}\\
\url{http://www.maroon.dti.ne.jp/fermat/a^4+hb^4=c^4+d^4/19000.pdf}\\

\end{document}